\documentclass{article}

\usepackage{graphicx}
\graphicspath{ {./images/} }

\begin{document}
\title{\centering \large \bf
New Flat surfaces in $S^3$}
\author{\large \bf
Armando M. V. Corro\\
IME, Universidade Federal de Goi\'as \\
Caixa Postal 131, 74001-970, Goi\^ania, GO, Brazil\\
e-mail:corro@ufg.br
  \and
 \large \bf Marcelo Lopes Ferro    \\
 IME, Universidade Federal de Goi\'as \\
Caixa Postal 131, 74001-970, Goi\^ania, GO, Brazil\\
e-mail: marceloferro@ufg.br
 }
\date{}
\maketitle \thispagestyle{empty}
\begin{abstract}
In this paper, we consider a method of constructing flat
surfaces based on Ribaucour transformations in the sphere 3-space. By applying the theory to the flat torus, we obtain a families of complete flat surfaces in $S^3$ which are determined by several parameters. we provide explicit examples.
\end{abstract}

\section*{\large \bf  Introduction}

Ribaucour transformations for hypersurfaces, parametrized by lines of curvature, were classically studied by Bianchi \cite{Bianchi}. They can be applied to obtain surfaces of constant Gaussian curvature and surfaces of constant mean curvature,  from a given such surface, respectively, with constant Gaussian curvature and constant mean curvature. The first application of this method to minimal and
cmc surfaces in $R^3$ was obtained by Corro, Ferreira, and Tenenblat in \cite{CFT1}-\cite{CFT3}. In \cite{TW}, Tenenblat and Wang extended such transformations to surfaces of space forms. For more application this method, see \cite{CT}, \cite{NT}, \cite{FT}, \cite{Tenenblat2}, \cite{TL} and \cite{TW1}.

Using Ribaucour transformations and applying the theory to the rotational flat surfaces in $H^3$, in \cite{CMT} the authours obtained families of new such surfaces.

The study of flat surfaces in $S^3$ traces back to Bianchi’s works in the 19th century, and it has a very rich global theory, as evidenced by the existence of a large class of flat tori in $S^3$,  see \cite{K1}, \cite{pinkall}, \cite{Spivak} and \cite{W}. Indeed, these flat tori constitute the only examples of compact surfaces of constant curvature in space forms that are not totally umbilical round spheres. Flat surfaces in $S^3$ admit a more explicit treatment than other constant curvature surfaces. Moreover, there are still important open
problems regarding flat surfaces in $S^3$, some of them unanswered for more than 40 years. For instance, it remains unknown if there exists an isometric embedding of R2 into $S^3$. These facts show that the geometry of flat surfaces in $S^3$ is a worth studying topic, although the number of contributions to the theory is not too large. Some 
important references of the theory are \cite{dadok}, \cite{K1}-\cite{K5}, \cite{pinkall}, \cite{Spivak} and \cite{W}.

In \cite{MJP}, the authors give a complete classification of helicoidal flat surfaces in $S^3$ by means of asymptotic lines coordinates.

In \cite{Aledo}, the authors characterized the flat surfaces in the unit 3-sphere that pass through a given regular curve of $S^3$ with a prescribed tangent plane distribution along this curve.

In this paper, motivated by \cite{CMT}  we use the  Ribaucour transformations to get a family of complete flat surfaces in $S^3$ from a given such surface in $S^3$. As an application of the theory, we obtain a families of complete flat surfaces in $S^3$ associated to the flat torus. The families we obtain depend four parameters. One of these parameters is given by parameterizing of the flat torus. The other parameters, appear from integrating the Ribaucour transformation. We show explicit examples of these surfaces. This work is organized as follows. In Section 1, we give a brief description of Ribaucour transformations in space forms. In Section 2, we give an additional condition for the transformed surface to be flat. In Section 3, we describe all flat surfaces of the sphere 3-space obtained by applying the
Ribaucour transformation to the flat torus. We prove that such surfaces are complete and provide explicit examples.
%\maketitle

\section*{ \large \bf 1. Preliminary}

This section contains the definitions and the basic theory of Ribaucour transformations for surfaces in $S^3$. ( For more Details see \cite{TW} )

Let $M$ be an orientable surface in $S^3$ without umbilic points, with Gauss map we denote by $N$. Suppose that 
there exist 2 orthonormal principal vector fields $e_1$ and $e_2$ defined on $M$. We say that $\widetilde{M}\subset S^3$ is associated to $M$ by a Ribaucour transformation with respect to $e_1$ and $e_2$, if there exist a differentiable function $h$ defined on $M$ and a
diffeomorphism $\psi:M\rightarrow \widetilde{M}$ such that\\
(a) for all $p\in M$,
$exp_ph(p)N(p)=exp_{\psi(p)}h(p)\widetilde{N}(\psi(p))$, where
$\widetilde{N}$ is the Gauss map of $\widetilde{M}$ and $exp$ is the expo-
nential map of $M^3(\overline{k})$.\\
(b)The subset $\{exp_ph(p)N(p)$, $p\in M\}$, is a two-dimensional
submanifold of $S^3$.\\
(c) $d\psi(e_i)$ $1\leq i\leq 2$ are orthogonal principal directions of $\widetilde{M}$.

\vspace{.1in}

The following result gives a characterization of Ribaucour
transfomations. ( see \cite{TW} for a proof and more details)

\vspace{.1in}

\noindent {\bf Theorem 1.1} \textit{ Let $M$ be an orientable
surface of $S^3$ parametrized by $X:U\subseteq R^2\rightarrow M$, without umbilic points. Assume $e_i=\frac{X,_i}{a_i}$, $1\leq i\leq 2$ where $a_i=\sqrt{g_{ii}}$ are orthogonal principal directions, $-\lambda_i$ the corresponding principal curvatures, and $N$ is a unit vector field normal to $M$. A surface $\widetilde{M}$ is locally associated to $M$ by a
Ribaucour transformation if and only if there is differentiable
functions $W,\Omega,\Omega_i:V\subseteq U\rightarrow R$ which
satisfy 
\begin{eqnarray}
\Omega_{i,j}&=&\Omega_j\frac{a_{j,i}}{a_{i}},\hspace{0,5cm}
for\hspace{0,2cm}i\neq j, \nonumber\\
\Omega,_i &=&a_i\Omega_i,\label{eq12}\\
W,_i&=&-a_i\Omega_i\lambda_i.\nonumber
\end{eqnarray}
$W(W+\lambda_i\Omega)\neq 0$ and $\widetilde{X}:V\subseteq
U\rightarrow \widetilde{M}$, is a parametrization of
$\widetilde{M}$ given by
\begin{eqnarray}
\widetilde{X}=\bigg(1-\frac{2\Omega^2}{S}\bigg)X-\frac{2\Omega}{S}\bigg(\sum_{i=1}^2\Omega_ie_i-WN\bigg),\label{eq6}
\end{eqnarray}
where
\begin{eqnarray}
S=\sum_{i=1}^2\big(\Omega_i\big)^2+W^2+\Omega^2.\label{eq7}
\end{eqnarray}
Moreover, the normal map of $\widetilde{X}$ is given by
\begin{eqnarray}
\widetilde{N}=N+\frac{2W}{S}\bigg(\sum_{i=1}^2\Omega_ie_i-WN+\Omega X\bigg),\label{eq8}
\end{eqnarray}
and the principal curvatures and coefficients of the first
fundamental form of $\widetilde{X}$, are given by
\begin{eqnarray}
\widetilde{\lambda}_i=\frac{WT_i+\lambda_iS}{S-\Omega
T_i},\hspace{1,5cm}\widetilde{g}_{ii}=\bigg(\frac{S-\Omega
T_i}{S}\bigg)^2g_{ii}\label{eq9}
\end{eqnarray}
where $\Omega_i$, $\Omega$ and $W$ satisfy (\ref{eq12}), $S$ is given by (\ref{eq7}), $g_{ii}$,
$1\leq i\leq 2$ are coefficients of the first fundamental form of
$X$, and
\begin{eqnarray}
T_1=2\bigg(\frac{\Omega_{1,1}}{a_1}+\frac{a_{1,2}}{a_1a_{2}}\Omega_2-W\lambda_1+\Omega\bigg),\hspace{0,1cm}T_2=2\bigg(\frac{\Omega_{2,2}}{a_2}+\frac{a_{2,1}}{a_{1}a_2}\Omega_1-W\lambda_2+\Omega\bigg)\label{T1T2first}
\end{eqnarray}
}

\section*{ \large \bf 2. Ribaucour transformation for flat surface in $S^3$}

In this section we provides a sufficient condition for a Ribaucour transformation to transform a flat surface into another such surface.

\vspace{.1in}

\noindent {\bf Theorem 2.1} \textit{ Let $M$ be a surfaces of $S^3$
parametrized by $X:U\subseteq R^2\rightarrow M$, without umbilic
points and let $\widetilde{M}$  parametrized by (\ref{eq6}) be
associated to $M$ by a Ribaucour transformation, such that the
normal lines intersect at a distance function $h$. Assume that
$h=\frac{\Omega}{W}$ is not constant along the lines of curvature
and the function $\Omega$, $\Omega_i$ and $W$ satisfy one of the
additional relation
\begin{eqnarray}
 \Omega_1^2+\Omega_2^2=c\big(\Omega^2+W^2\big)\label{eq13}
\end{eqnarray}
where $c>0$, $S$ is given by (\ref{eq7}) and $W$, $\Omega_i$, $1\leq i\leq 2$ satisfies (\ref{eq12}). Then $\widetilde{M}$ parameterized by
(\ref{eq6}) is a flat surface, if and only if $M$ is a
flat surface.
 }

\vspace{.1in}

\noindent {\bf Proof:} See \cite{TW} Theorem 2.1, with $\alpha=\gamma=1$ and $\beta=0$.
%\mbox{}\hfill $\Box$

\vspace{.1in}

 \noindent {\bf Remark 2.2}  Let $X$ as in the previous theorem.
Then the parameterization $\widetilde{X}$ of $\widetilde{M}$,
locally associated to $X$ by a Ribaucour transformation, given by
(\ref{eq6}), is defined on
\begin{eqnarray*}
V=\{(u_1,u_2)\in U;\hspace{0,1cm}(\Omega T_1-S)(\Omega T_2-S)\neq0\}.\nonumber
\end{eqnarray*}

\vspace{.1in}

\section*{ \large \bf 3. Families of flat surfaces associated to the Flat Torus in $S^3$.}

\vspace{.2in}

In this section, by applying Theorem 2.1 to the Flat Torus in $S^3$, we obtain a three parameter family of complete flat surfaces in $S^3$.

 \vspace{.2in}

 \noindent {\bf Theorem 3.1} \textit{ Consider the Flat Torus in $S^3$ parametrized by
\begin{eqnarray}
X(u_1,u_2)=(r_1\cos(r_2u_1),r_1\sin(r_2u_1),r_2\cos(r_1u_2),r_2\sin(r_1u_2)), \hspace{0,3cm}(u_1,u_2)\in
R^2\label{torus}
\end{eqnarray}
$r_i$, $1\leq i\leq 2$ are positive constants satisfying $r_1^2+r_2^2=1$,
as flat torus in $S^3$ where the first fundamental form is
$I=r_1^2r_2^2\big(du_1^2+du_2^2\big)$. A parametrized surface $\widetilde{X}(u_1,u_2)$
is flat surface locally associated to $X$ by a Ribaucour
transformation as in Theorem 2.1 with additional relation given by (\ref{eq13}), if and only if, up to an isometries of $S^3$, it is given by
\begin{eqnarray}\label{xtiu}
\widetilde{X}=\frac{1}{S}\left[  \begin{array}{c}
   \big(r_1S-2r_1\Omega^2-2r_2\Omega W\big)\cos(r_2u_1)+2f'\Omega\sin(r_2u_1)\\
    \big(r_1S-2r_1\Omega^2-2r_2\Omega W\big)\sin(r_2u_1)-2f'\Omega\cos(r_2u_1)\\
    \big(r_2S-2r_2\Omega^2+2r_1\Omega W\big)\cos(r_1u_2)+2g'\Omega\sin(r_1u_2)\\
    \big(r_2S-2r_2\Omega^2+2r_1\Omega W\big)\sin(r_1u_2)-2g'\Omega\cos(r_1u_2)\\
    \end{array} \right]
\end{eqnarray}
defined on $V=\{(u_1,u_2)\in R^2;\hspace{0,1cm}(r_1^2g^2-r_2^2f^2-2r_2^2fg)(r_2^2f^2-r_1^2g^2+2r_1^2fg)\neq0\}$, where
\begin{eqnarray}
\Omega=r_1r_2\big(f(u_1)+g(u_2)\big),\,\,\, W=r_2^2f(u_1)-r_1^2g(u_2),\,\,\,S=(1+c)\big(\Omega^2+W^2\big)\label{omegaeW}
\end{eqnarray}
$c> 0$, and the functions $f$ and $g$ are given by
\begin{eqnarray}
&&i)\,\,f(u_1)=\cosh(r_2\sqrt{c}u_1),\,\,\,g(u_2)=\frac{r_2}{r_1}\sinh(r_1\sqrt{c}u_2),\label{fgA1positivo}\\
&& \hspace{4cm}or\nonumber\\
&&ii)\,\,f(u_1)=\sinh(r_2\sqrt{c}u_1),\,\,\,g(u_2)=\frac{r_2}{r_1}\cosh(r_1\sqrt{c}u_2),\label{fgA1negativo}\\\
&&\hspace{4cm} or\nonumber\\
&&iii)\,\,f(u_1)=a_1e^{\epsilon_1r_2\sqrt{c}u_1},\,\,\,g(u_2)=b_1e^{\epsilon_2r_1\sqrt{c}u_2},\,\, \epsilon_i^2=1, \, 1\leq i\leq 2.\label{fgA1zero}
\end{eqnarray}
Moreover, the normal map of $\widetilde{X}$ is given by
\begin{eqnarray}\label{ntiu}
\widetilde{N}=\frac{1}{S}\left[  \begin{array}{c}
   \big(-r_2S-2r_2W^2+2r_1\Omega W\big)\cos(r_2u_1)+2f'\Omega\sin(r_2u_1)\\
    \big(-r_2S-2r_2W^2+2r_1\Omega W\big)\sin(r_2u_1)-2f'\Omega\cos(r_2u_1)\\
    \big(r_1S+2r_1W^2-2r_2\Omega W\big)\cos(r_1u_2)+2g'\Omega\sin(r_1u_2)\\
    \big(r_1S+2r_1W^2-2r_2\Omega W\big)\sin(r_1u_2)-2g'\Omega\cos(r_1u_2)\\
    \end{array} \right]
\end{eqnarray}
}

\vspace{.1in}

\noindent {\bf Proof:} Consider the first fundamental form of the
Flat Torus $ds^2=r_1^2r_2^2(du_1^2+du_2^2)$ and the principal curvatures $-\lambda_i$ $1\leq i\leq 2$ given by
$\lambda_1=\frac{-r_2}{r_1}$, $\lambda_2=\frac{r_1}{r_2}$. Using (\ref{eq12}), to obtain the
Ribaucour transformations, we need to solve the following of
equations
\begin{eqnarray}
\Omega_{i,j}=0,\hspace{0,5cm}\Omega,_i
=r_1r_2\Omega_i,\hspace{0,5cm}W,_i=-r_1r_2\Omega_i\lambda_i,
\hspace{0,2cm}1\leq i\neq j\leq 2.\label{eqomega}
\end{eqnarray}
Therefore we obtain
\begin{eqnarray}
\Omega =r_1r_2\big(f_1(u_1)+f_2(u_2)\big),\hspace{0,5cm}W=-r_1r_2\big(\lambda_1f_1+\lambda_2f_2\big)+\overline{c},\hspace{0,5cm}\Omega_i=f_i',
\hspace{0,2cm}1\leq i\neq j\leq 2,\label{omega}
\end{eqnarray}
 where $\overline{c}$ is a real constant. Using (\ref{eq13})
the associated surface will be flat
when $\displaystyle{\Omega_1^2+\Omega_2^2=c(\Omega^2+W^2)}$. Therefore, we obtain
that $c> 0$ and the functions $f_1$ and $f_2$ satisfy
\begin{eqnarray}
(f_1')^2+(f_2')^2=c(\Omega^2+W^2).\label{eqf1''}
\end{eqnarray}
Differentiate this last equation with respect $x_1$ and $x_2$, using (\ref{eqomega}) and (\ref{omega}) we get
\begin{eqnarray}
&&f_1''=cr_2^2f_1+cr_2^2\overline{c},\nonumber\\
&&f_2''=cr_1^2f_2-cr_1^2\overline{c}.\nonumber
\end{eqnarray}
Defining $f(u_1)=f_1(u_1)+\overline{c}$ and $g(u_2)=f_2(u_2)-\overline{c}$, we have
\begin{eqnarray}
f''-cr_2^2f=0,&&\,\,\,\,\,\ g''-cr_1^2g=0,\label{eqf1'''}\\
\Omega =r_1r_2\big(f(u_1)+g(u_2)\big),&&\,\,\,\,\,\ W=r_2^2f(u_1)-r_1^2g(u_2).\label{omegafim}
\end{eqnarray}
By Theorem 1.1, we have that $\widetilde{X}$ and $\widetilde{N}$ are
given by (\ref{xtiu}) and (\ref{ntiu}). Using (\ref{T1T2first}) and (\ref{omegafim}), we get
$$T_1=\frac{2r_2(1+c)f}{r_1}\,\,\,\,\,\,\,T_2=\frac{2r_1(1+c)g}{r_2}.$$
Thus, from Remark 2.2 $\widetilde{X}$ is defined
in $V=\{(u_1,u_2)\in R^2;\hspace{0,1cm}(r_1^2g^2-r_2^2f^2-2r_2^2fg)(r_2^2f^2-r_1^2g^2+2r_1^2fg)\neq0\}$.\\

From (\ref{eqf1'''}), we get
\begin{eqnarray}
f(u_1)=a_1\cosh(r_2\sqrt{c}u_1)+a_2\sinh(r_2\sqrt{c}u_1),\,\,\ g(u_2)=b_1\cosh(r_1\sqrt{c}u_2)+b_2\sinh(r_1\sqrt{c}u_2). \label{fgprimeiro}
\end{eqnarray}
Substituting (\ref{fgprimeiro}) and (\ref{omegafim}) in (\ref{eqf1''}), we have
\begin{eqnarray}
\big(a_1^2-a_2^2\big)r_2^2=\big(b_2^2-b_1^2\big)r_1^2.\label{eqa1b1}
\end{eqnarray}
Let $A_1=a_1^2-a_2^2$. If $A_1>0$, then from (\ref{eqa1b1}) $b_2^2> b_1^2$. Hence (\ref{fgprimeiro}) can be rewritten as
\begin{eqnarray}
&&f(u_1)=\sqrt{A_1}\cosh(r_2\sqrt{c}u_1+A_2),\label{fterceiro}\\
&&g(u_2)=\sqrt{A_1}\frac{r_2}{r_1}\sinh(r_1\sqrt{c}u_2+B_2),\label{gterceiro}
\end{eqnarray}
where
$\cosh(A_2)=\frac{a_1}{\sqrt{a_1^2-a_2^2}}$, $\sinh(A_2)=\frac{a_2}{\sqrt{a_1^2-a_2^2}}$, $\sinh(B_2)=\frac{b_2}{\sqrt{b_2^2-b_1^2}}$ and $\cosh(B_2)=\frac{b_1}{\sqrt{b_2^2-b_1^2}}$.\\

The constants $A_2$ and $B_2$, without loss of
generality, my be considered to be zero. One can verify that the surfaces with different values of $A_2$ and $B_2$
are congruent. In fact, using the notation $\widetilde{X}_{A2B2}$
for the surface
$\widetilde{X}$ with fixed constants $A_2$ and $B_2$, we have
$$\widetilde{X}_{A2B2}=R_{(\frac{-A_2}{r_2\sqrt{c}},\frac{-B_2}{r_1\sqrt{c}})}\widetilde{X}_{00}\circ h $$
where $\displaystyle{h(u_1,u_2)=\bigg(u_1+\frac{A_2}{r_2\sqrt{c}},u_1+\frac{B_2}{r_1\sqrt{c}}\bigg)}$
with 
$$R_{(\theta,\phi)}(x_1,x_2,x_3,x_4)=(x_1\cos\theta-x_2\sin\theta,x_1\sin\theta+x_2\cos\theta,x_3\cos\phi-x_4\sin\phi,x_3\sin\phi+x_4\cos\phi).$$
Now substituting (\ref{fterceiro}) with $A_2=0$, (\ref{gterceiro}) with $B_2=0$ and (\ref{omegafim}) in (\ref{xtiu}) we obtain that $\widetilde{X}$ does not depend on $A_1$. Thus without loss of generality, we can consider $A_1=1$. Therefore we conclude that $f$ and $g$ are given by (\ref{fgA1positivo}).

On the other hand, if $A_1<0$, then from (\ref{eqa1b1}) $b_2^2< b_1^2$. Hence (\ref{fgprimeiro}) can be rewritten as
\begin{eqnarray}
&&f(u_1)=\sqrt{-A_1}\sinh(r_2\sqrt{c}u_1+A_2),\label{fterceiro1}\\
&&g(u_2)=\sqrt{-A_1}\frac{r_2}{r_1}\cosh(r_1\sqrt{c}u_2+B_2).\label{gterceiro1}
\end{eqnarray}
Proceeding in a similar way to the previous case, we obtain that $f$ and $g$ are given by (\ref{fgA1negativo}).

If $A_1=0$, then $a_2=\epsilon_1a_1$, and from (\ref{eqa1b1}) $b_2=\epsilon_2b_1$, with $\epsilon_i^2=1$, $1\leq i\leq 2$. Thus, substituting this in (\ref{fgprimeiro}), we obtain (\ref{fgA1zero}).

\vspace{.2in}

 \noindent {\bf Remark 3.2} Each flat surfaces associated to
 the flat torus as in Theorem 3.1, is parametrized by lines of
 curvature and from (\ref{eq9}), the metric is given by
 $ds^2=\psi_1^2du_1^2+\psi_2^2du_2^2$, where
\begin{eqnarray}
\psi_1=\frac{(-r_2^2f^2+r_1^2g^2-2r_2^2fg)r_1r_2}{r_1^2g^2+r_2^2f^2},\,\,\,\psi_2=\frac{(r_2^2f^2-r_1^2g^2-2r_1^2fg)r_1r_2}{r_1^2g^2+r_2^2f^2}.\label{metrica}
\end{eqnarray}
Moreover, from (\ref{eq9}), the principal curvatures of the
$\widetilde{X}$ are given by
\begin{eqnarray}
\widetilde{\lambda}_1=\frac{-r_2\psi_2}{r_1\psi_1},\,\,\,\,\,\,\widetilde{\lambda}_2=\frac{r_1\psi_1}{r_2\psi_2}.\label{aa1}
\end{eqnarray}

\vspace{.1in}

\noindent{\bf Proposition 3.3} \textit{Any flat surfaces
associated to the flat torus $\widetilde{X}$, given by
Theorem 3.1 is complete.}

\vspace{.15in}

\noindent {\bf Proof:}
For divergent curves
$\gamma(t)=(u_1(t),u_2(t))$, such that
$\displaystyle{\lim_{t\rightarrow
\infty}\big(u_1^2+u_2^2\big)=\infty}$, we have
$l(\widetilde{X}\circ\gamma)=\infty$.

In fact, the functions $f$ and $g$ are given by
(\ref{fgA1positivo}) or (\ref{fgA1negativo}) or (\ref{fgA1zero}) and the coefficients of the first fundamental form  $\psi_i$, $1\leq i\leq 2$ of $\widetilde{X}$ given by (\ref{metrica}).
Therefore $\displaystyle{\lim_{|u_1|\rightarrow
\infty}|\psi_i|=r_1r_2}$, $1\leq i\leq 2$, uniformly in $u_2$ and  $\displaystyle{\lim_{|u_2|\rightarrow
\infty}|\psi_i|=r_1r_2}$, $1\leq i\leq 2$, uniformly in $u_1$. Hence, there are $k_1>0$ and $k_2>0$ such that $|\psi_i(u_1,u_2)|>\frac{r_1r_2}{2}$, $1\leq i\leq 2$  for all $(u_1,u_2)\in R^2$ with $|u_1|>k_1$ and $|u_2|>k_2$.\\
Let
$$m_i=min\{|\psi_i(u_1,u_2)|;\,(u_1,u_2)\in [-k_1,k_1 ]\times [-k_2,k_2 ]\}.$$
Therefore $|\psi_i(u_1,u_2)|\geq m_i$ in $[-k_1,k_1 ]\times [-k_2,k_2 ]$. Now consider $m_0=min\{m_1,m_2,\frac{r_1r_2}{2}\}$, then  $|\psi_i(u_1,u_2)|\geq m_0$ in $R^2$. We conclude that $\widetilde{X}$ is a complete surface.

 \vspace{.2in}

In the following, we provide some examples

\vspace{.2in}

\noindent {\bf Example 3.3}{ Consider stereographic projection $\pi:S^3\rightarrow R^3$
$$\pi(x_1,x_2,x_3,x_4)=\frac{1}{1-x_4}\bigg(x_1,x_2,x_3\bigg)$$
where $x_1^2+x_2^2+x_3^2+x_4^2=1$.\\

In Figure 1, we provide the surface parametrized by $\pi\circ\widetilde{X}$ where $\widetilde{X}$ is given by (\ref{xtiu}) locally associated to flat torus in $S^3$ by a Ribaucour transformation. In this surface, we have\\
$f(u_1)=\cosh\bigg(\frac{8u_1}{5}\bigg)$ and $g(u_2)=\frac{4}{3}\sinh\bigg(\frac{6u_2}{5}\bigg)$\\
\begin{figure}[!h]
\begin{center}
{\includegraphics[scale=0.5]{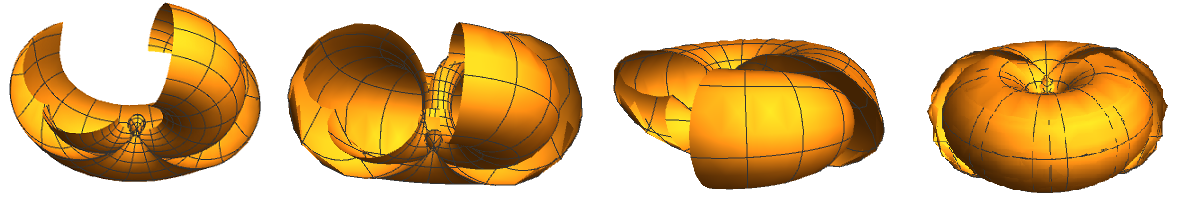}}
\caption[!h]{In the figure above we have $r_1=\frac{3}{5}$, $r_2=\frac{4}{5}$, $c=4$.}
\end{center}
\end{figure}}

\vspace{.1in}

In Figure 2, we provide the surface parametrized by $\pi\circ\widetilde{X}$ where $\widetilde{X}$ is given by (\ref{xtiu}) locally associated to flat torus in $S^3$ by a Ribaucour transformation. In this surface, we have\\
$f(u_1)=\sinh\bigg(\frac{8u_1}{5}\bigg)$ and $g(u_2)=\frac{4}{3}\cosh\bigg(\frac{6u_2}{5}\bigg).$\\
\begin{figure}[!h]
\begin{center}
{\includegraphics[scale=0.5]{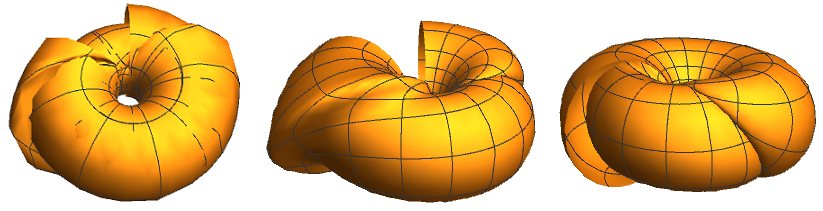}}
\caption[!h]{In the figure above we have $r_1=\frac{3}{5}$, $r_2=\frac{4}{5}$, $c=4$.}
\end{center}
\end{figure}

\vspace{.1in}

In Figure 3, we provide two surfaces parametrized by $\pi\circ\widetilde{X}$ where $\widetilde{X}$ is given by (\ref{xtiu}) locally associated to flat torus in $S^3$ by a Ribaucour transformation. In this surfaces, we have $r_1=\frac{3}{5}$, $r_2=\frac{4}{5}$, $c=\frac{1}{1000}$ with $f$ and $g$ given by (\ref{fgA1positivo})  in the first surface
and $f$ and $g$ given by (\ref{fgA1negativo})  and the second.
\begin{figure}[!h]
\begin{center}
{\includegraphics[scale=0.5]{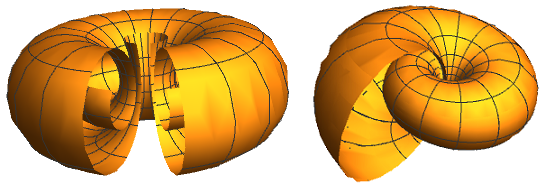}}
\caption[!h]{}
\end{center}
\end{figure}

\vspace{.1in}

In Figure 4, we provide two surfaces parametrized by $\pi\circ\widetilde{X}$ where $\widetilde{X}$ is given by (\ref{xtiu}) locally associated to flat torus in $S^3$ by a Ribaucour transformation. In this surfaces, we have $r_1=\frac{3}{5}$, $r_2=\frac{4}{5}$, $c=\frac{1}{1000}$, $a_1=b_1=1$ with $f$ and $g$ given by (\ref{fgA1zero}) .
\begin{figure}[!h]
\begin{center}
{\includegraphics[scale=0.5]{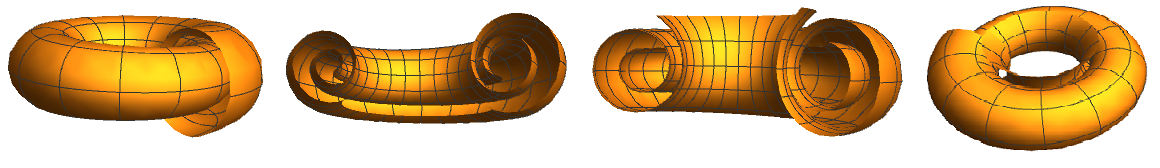}}
\caption[!h]{}
\end{center}
\end{figure}

{}

\end{document}